\documentclass[12pt,oneside]{amsart}
\usepackage{times}
\usepackage{chicago}

\input BoxedEPS.tex
\SetepsfEPSFSpecial
\HideDisplacementBoxes

\catcode`\@=11

\def\cases@#1{\left\{\,\vcenter{\normalbaselines\m@th
    \ialign{$##$\hfil&\quad##\hfil\crcr#1\crcr}}\right.}

\def\Cases#1\endCases{\let\\=\cr\cases@{#1}}
\catcode`@=12

\def\subheading#1{\par\noindent{\bf #1.}}
\def\demo#1{\subheading{#1}}

\theoremstyle{plain}
\newtheorem{theorem}{Theorem}
\newtheorem{lemma}[theorem]{Lemma}
\newtheorem{definition}[theorem]{Definition}
\newtheorem{corollary}[theorem]{Corollary}

\theoremstyle{definition}
\newtheorem{example}[theorem]{example}

\theoremstyle{remark}
\newenvironment{remark}{\begin{quote}{\bf Remark.\quad}\rm\small}{\end{quote}}

\theoremstyle{remarks}
\newenvironment{remarks}{\begin{quote}{\bf Remarks.\quad}\rm\small}{\end{quote}}

\def\Corollary{\begin{corollary}}
\def\endCorollary{\end{corollary}}

\def\Definition{\begin{definition}}
\def\endDefinition{\end{definition}}

\def\Theorem{\begin{theorem}}
\def\endTheorem{\end{theorem}}

\def\NamedTheorem{\begin{theorem}}
\def\endNamedTheorem{\end{theorem}}

\def\Lemma{\begin{lemma}}
\def\endLemma{\end{lemma}}

\def\Example{\begin{example}}
\def\endExample{\end{example}}

\def\Remark{\begin{remark}}
\def\endRemark{\end{remark}}

\def\Remarks{\begin{remarks}}
\def\endRemarks{\end{remarks}}

\def\Proof{\begin{proof}}
\def\endProof{\end{proof}}

\let\Cal=\mathcal
\def\aa{{\Cal A}}

\def\pp{{\Cal P}}

\def\ww{{\Cal W}}\def\xx{{\Cal X}}

\def\RR{\mathbb{R}}
\def\PP{\mathbb{P}}
\def\QQ{\mathbb{Q}}
\def\text#1{\hbox{\rm#1}}
\def\qt#1{\qquad\text{#1}}

\def\siln{{\textstyle\sum\nolimits_{i\le n}}}
\def\NL{\nolimits}

\def\norm#1{\left\vert#1\right\vert}
\def\cref#1{(\ref{#1})}
\let\sref=\ref
\def\wrt/{with respect to}
\def\ritem(#1){\item}

\def\align{\begin{array}{rl}}
\def\endalign{\end{array}}

\def\th{\theta}
\def\gam{\gamma}
\def\del{\delta}

\def\thetabar{\bar\theta}

\def\siln{\sum\NL_{i\le n}}
\def\co#1{\text{co}\left(#1\right)}

\def\psiprop/{{\smc good}-$\psi$}

\begin{document}

\title[Talagrand's concentration inequality]{A note on Talagrand's convex hull
concentration inequality}

\author{David Pollard}
\address{
Statistics Department\\
 Yale University\\
 Box 208290 Yale Station\\
 New Haven, CT 06520-8290.}
\email{david.pollard@yale.edu}
\urladdr{http://www.stat.yale.edu/\~{}pollard/}


\date{28 September 2005}


\begin{abstract}\normalsize
The paper reexamines an argument by Talagrand that leads to a remarkable exponential
tail bound for the concentration of probability near a set.  The main novelty is the
replacement of a mysterious Calculus inequality by an application of Jensen's inequality.

\end{abstract}
\maketitle

\section{Introduction}  \label{intro}
Let $\xx$ be a set equipped with a sigma-field~$\aa$.  For each vector $w=(w_1,\dots,w_n)$
in~$\RR_+^n$, the weighted Hamming distance between two vectors
$x=(x_1,\dots,x_n)$ and $y=(y_1,\dots,y_n)$, in~$\xx^n$ is defined as
$$
d_w(x,y) := \sum_{i\le n}w_ih_i(x,y)\qt{where }
h_i(x,y)
=\begin{cases}
$1$&\text{if }x_i\ne y_i\cr
$0$&\text{otherwise}.
\end{cases}
$$
For a subset~$A$ of~$\xx^n$ and $x\in\xx^n$, the distances $d_w(x,A)$ and $D(x,A)$ are
defined by
$$
d_w(x) := \inf\{y\in A: d_w(x,y)\}
$$
and
$$
D(x,A) := \sup\NL_{w\in\ww} d_w(X,A),
$$
where the supremum is taken over all weights in the set
$$
\ww := \{(w_1,\dots,w_n): w_i\ge0\text{ for each $i$  and }\norm{w}^2 := \siln w_i^2\le1\}.
$$

\citeN[Section~4.1]{Talagrand95IHES} proved a remarkable concentration inequality for
 random elements~$X=(X_1,\dots,X_n)$ of~$\xx^n$ with independent
coordinates and  subsets~$A\in\aa^n$:
\begin{equation}
\PP\{X\in A\}\PP\{D(X,A)\ge t\}\le\exp(-t^2/4)\qt{for all
$t\ge0$}.
\label{Talagrand.4.1.1}
\end{equation}
As Talagrand showed, this inequality has many applications to problems in combinatorial
optimization and other other areas. See  \citeN{Talagrand96newlook},
\citeN[Chapter~6]{Steele97siam}, and
\citeN[Section~4]{McDiarmid98conc} for further examples.

Talagrand used an induction on~$n$ to establish his result, invoking a slightly
mysterious Calculus  lemma in the inductive step.  
%Apparently
%some readers regard  inductive arguments as less satisfactory than  direct proofs.
There has been a strong push in the literature to establish concentration and deviation
inequalities by ``more intuitive'' methods, such as those based on the tensorization, as in
\citeN{Ledoux96ESAIM},
\citeN{BoucheronLugosimassart2000ras},
\citeN{Massart03Flour}, and \citeN{Lugosi2003ANU}.

It is my purpose in this note to modify Talagrand's proof---adapting an idea from
\citeN[Section~3]{Talagrand96IM}---so that the inductive  step becomes a simple
application of the H\"older inequality (essentially as in the original proof) and the Jensen
inequality.

The distance
$D(x,A) $ has another representation, as a minimization over a convex subset of~$[0,1]$. 
Write~$h(x,y)$ for the point of~$\{0,1\}^n$ with~$i$th coordinate~$h_i(x,y)$.  For each
fixed~$x$, the function
$ h(x,\cdot)$  maps
$A$ onto a subset
$
h(x,A) := \{h(x,y):y\in A\}
$
of~$\{0,1\}^n$. The convex hull~$\co{h(x,A)}$ of~$h(x,A)$ in~$[0,1]^n$ is compact, and
$$
D(x,A) = \inf\{|\xi|: \xi\in \co{h(x,A)}\}  .
$$
Each point~$\xi$  of~$\co{h(x,A)}$ can be written as  $\int
h(x,y)\,\nu(dy)$ for a~$\nu$ in the set~$\pp(A)$ of all Borel probability measures for
which~$\nu(A)=1$. That is, $\xi_i = \nu\{y\in A: y_i\ne x_i\}$.  Thus
\begin{equation}
D(x,A)^2 =
\inf_{\nu\in \pp(A)} \siln \left(\strut\nu\{y\in A: y_i\ne x_i\}\right)^2  .
\label{DxA2}
\end{equation}

Talagrand actually proved inequality~\cref{Talagrand.4.1.1}
by showing that
\begin{equation}
\PP \{X\in A\}\PP \exp\left(\strut \tfrac14 D(X,A)^2\right) \le 1.
\label{expD2}
\end{equation}
He also established  an even stronger result, in which
the $D(X,A)^2/4$ in~\cref{expD2} is replaced by a more complicated distance
function. 

For each convex, increasing function~$\psi$ with $\psi(0)=0=\psi'(0)$ define
\begin{equation}
F_\psi(x,A) := \inf_{\nu\in\pp(A)}\siln \psi\left(\strut \nu\{y\in A: y_i\ne x_i\}\right),
\label{Fpsi}
\end{equation}
%The inequality~\cref{expD2} corresponds to the choice $\psi(\th)=\th^2/4$.
For each $c>0$, \citeN[Section~4.2]{Talagrand95IHES} showed that
\begin{equation}
(\PP \{X\in A\})^c\PP \exp\left(\strut F_{\psi_c}(X,A)\right) \le 1,
\label{Fpsic}
\end{equation}
where
\begin{equation}\label{psic}
\begin{align}  
 \psi_c(\th) &:=
 c^{-1}\left((1-\th)\log (1-\th) - (1-\th+c)\log\left(\strut \dfrac{(1-\th)+c}{
1+c}\right)\right) \\
&=  \sum\NL_{k\ge 2}\dfrac{\th^k}{k}\left(\strut \dfrac{R_c +R_c^2+\dots +R_c^{k-1}}{
(k-1)}
\right)
\qt{with $R_c := \dfrac{1}{c+1}$.}\\
&\ge \dfrac{\th^2}{2+2c}
\end{align}
\end{equation}
As you will  see in Section~\ref{concave}, this strange function is actually the
largest solution to a differential inequality,
$$
\psi''(1-\th)\le {1/ (\th^2 +\th c)}\qt{for $0<\th<1$}.  %%\label{diff.ineq}
$$
 Inequality~\cref{Fpsic} improves
on~\cref{expD2} because
$
 D(x,A)^2 /4\le F_{\psi_1}(x,A)
$.

Following the lead of
\citeN[Section~4.4]{Talagrand95IHES}, we can ask for general conditions on the convex~$\psi$
under which an analog of~\cref{Fpsic} holds with some other decreasing function of $\PP
\{X\in A\}$ as an upper bound. The following modification of Talagrand's theorems gives a
sufficient condition in a form that serves to emphasize the role played by Jensen's
inequlity.

\Theorem \label{convex.conc}
Suppose $\gam$ is a decreasing function with $\gam(0)=\infty$ and $\psi$ is a
convex function.  Define $G(\eta,\th) := \psi(1-\th)+\th\eta$
and $G(\eta):=\inf_{0\le\th\le1}G(\eta,\th)$ for $\eta\in\RR^+$. 
Suppose
\begin{enumerate}
\item
$
r\mapsto \exp\left(\strut
G_\psi(\gam(r)-\gam(r_0))\right)$
 is  concave on~$[0,r_0]$, for each $r_0\le1$

\item
$
(1-p)e^{\psi(1)} + p \le e^{\gam(p)}$ for $0\le p\le 1$.
\end{enumerate}
Then
\begin{equation*}
\PP \exp\left(\strut F_\psi(X,A)\right) \le \exp\left(\strut \gam\left(\strut \PP \{X\in
A\}\right)\right).
\end{equation*}
for every $A\in\aa^n$ and every random element~$X$ of~$\xx^n$ with independent components.
\endTheorem

The following lemma, a more general version of which is proved 
 in Section~\sref{concave}, leads to a simple sufficient condition for the concavity
assumption~(ii)  of Theorem~\ref{convex.conc} to hold.

\begin{lemma}[Concavity lemma]  \label{concavity.lemma}
Suppose $\psi:[0,1]\to\RR^+$  is convex and increasing, with
$\psi(0)=0=\psi'(0)$ and $\psi''(\th)>0$ for $0<\th<1$. 
Suppose  $\xi:[0,r_0]\to\RR^+\cup\{\infty\}$ is continuous  and twice
differentiable on~$(0,r_0)$. Suppose also that there exists some finite constant~$c$ for
which
$\xi''(r)\le c\xi'(r)^2$ for $0<r<r_0$.  If
$$
\psi''(1-\th)\le {1/ (\th^2 +\th c)}\qt{for $0<\th<1$}
$$
then the function $r\mapsto \exp\left(\strut G(\xi(r))\right)$ is concave on~$[0,r_0]$.
\end{lemma}

The Lemma will be applied with $\xi(r)=\gam(r)-\gam(r_0)$ for $0\le r\le r_0$.  As shown 
in Section~\sref{concave}, the conditions of the Lemma hold for $\psi(\th)=\th^2/4$ with
$\gam(r)= \log(1/r)$ and also for the $\psi_c$ from~\cref{psic} with $\gam(r) = c\log(1/r)$.

\begin{remarks}
\begin{enumerate}
\item
If $\gam(0)$ were finite, the  inequality asserted by Theorem~\ref{convex.conc} could not
hold for all nonempty~$A$ and all~$X$. For example, if each~$X_i$ had a nonatomic
distribution and~$A$ were a singleton set we would have
$F_\psi(X,A)=n \psi(1)$ almost surely. The quantity $\PP \exp\left(\strut
F_\psi(X,A)\right) $ would exceed
$\exp(\gam(0))$ for large enough~$n$. It it to avoid this difficulty that we
need~$\gam(0)=\infty$.

\item
Assumption~(ii) of the Theorem, which is essentially an assumption that the asserted
inequality holds for~$n=1$, is easy to check if $\gam$ is a convex function with
$\gam(1)\ge0$. For then the  function $B(p):=\exp(\gam(p))$ is convex
 with $B(1)\ge1$ and 
$B'(1)= \gam'(1)e^{\gam(1)}$.  We have 
$$
B(p)\ge (1-p)e^{\psi(1)} +p 
\qt{for all $p$ in~$[0,1]$}
$$
 if $B'(1)\le 1-e^{\psi(1)}$.

\item
I had hoped to extend the proof to cover the case $c=0$ but I then ran into problems
with
$\gam(0)=\infty$.
\end{enumerate}\end{remarks}

\section{Proof of  Theorem~\ref{convex.conc}}  \label{GCHTproof}
Argue by induction on~$n$. 
As a way of keeping the notation straight, replace the subscript on~$F_\psi(x,B)$ by an~$n$
when the argument~$B$ is a subset of~$\xx^n$. Also, work with the product measure~$\QQ
=\otimes_{i\le n} Q_i$ for the distribution of~$X$ and $\QQ_{-n}
=\otimes_{i< n} Q_i$ for the distribution of $(X_1,\dots,X_{n-1})$. The assertion of the
Theorem then becomes
$$
\QQ \exp\left(\strut F_n(x,A)\right) \le \exp(\gam(\QQ A))
$$

For $n=1$ and~$B\in\aa$  we have
$
F_1(x,B) =  \psi(1)\{x\notin B\}+ 0\{x\in B\}
$
so that $Q_1\exp\left(\strut F_1(x,B)\right) \le (1-p)e^{\psi(1)} +p$, where $p=Q_1 B$.
Assumption~(i) then gives the desired~$\exp(\gam(p))$ bound.

Now suppose that $n>1$ and that the inductive hypothesis is valid for dimensions strictly
smaller than~$n$.
Write
$\QQ$ as
$\QQ_{-n}\otimes Q_n$. To simplify notation, write $w$ for
$x_{-n}:=(x_1,\dots,x_{n-1})$ and $z$ for~$x_n$. Define the cross section~
$
A_z:=\{w\in\xx^{n-1}:
(w,z)\in A\}
$
and write $R_z$ for $\QQ_{-n} A_z$. Define $r_0
:=\sup_{z\in\xx}R_z$.  Notice that 
$
r_0 \ge Q_n^z R_z = \QQ A
$.

The key to the proof is a recursive bound for~$F_n$: for each $x=(w,z)$
with~$A_z\ne\emptyset$, each~$m$ with $A_m\neq \emptyset$, and all  $\th\in[0,1]$,
\begin{equation} \label{recursive.F}
F_n(x,A) \le \th F_{n-1}(w,A_z) +\thetabar F_{n-1}(w,A_m) +\psi(\thetabar)
\qt{where }\thetabar:=1-\th .
\end{equation}

\begin{figure}[htb]
\centerline{
\BoxedEPSF{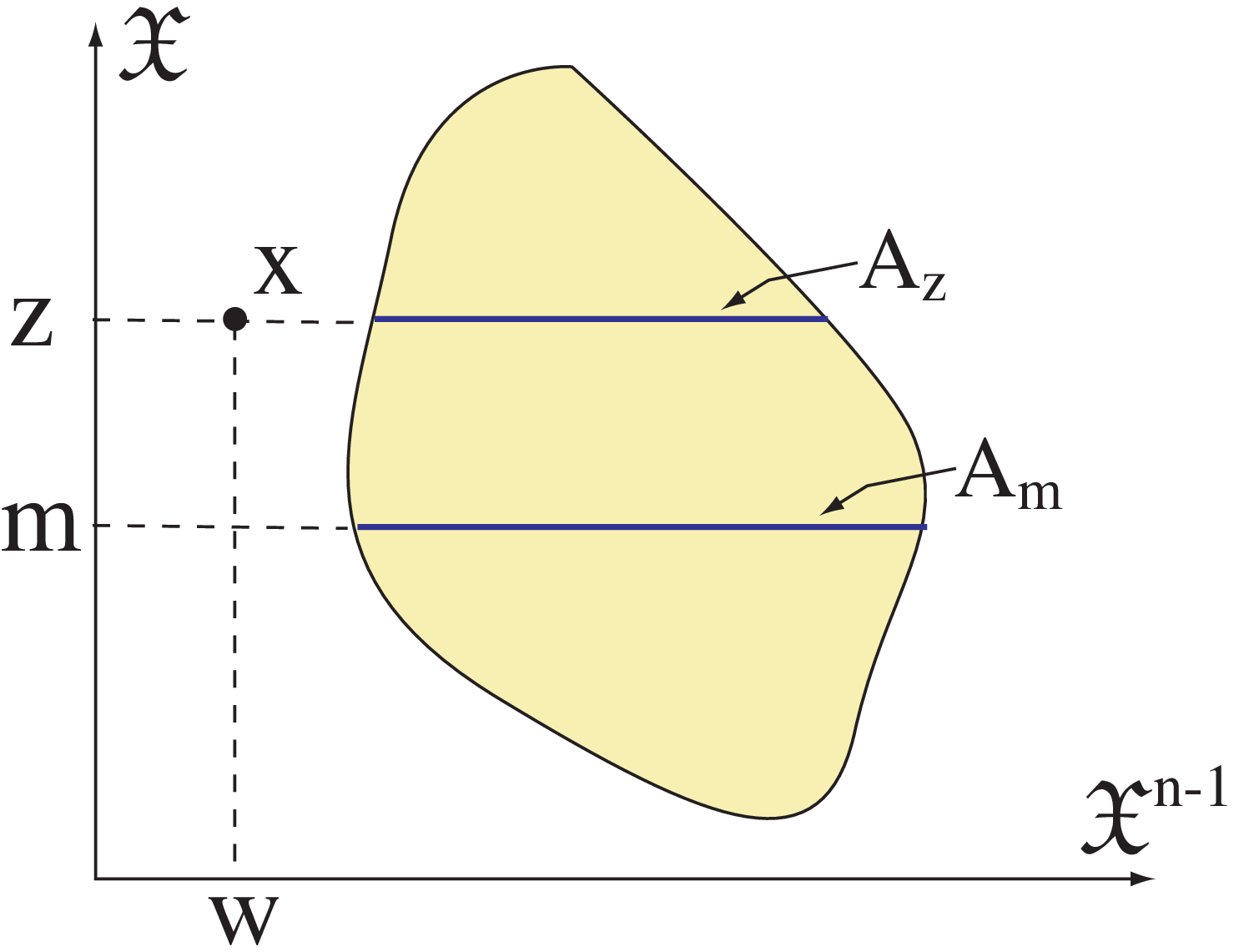 scaled 400}}
\end{figure}

To establish   inequality~\cref{recursive.F},
suppose  $\mu_z$ is a probability measure concentrated on~$A_z$ and $\mu_m$ is a probability
measure concentrated on~$A_m$.  
For a~$\th$ in~$[0,1]$, define
$
\nu = \th\mu_z\otimes\del_z + \thetabar\mu_m\otimes\del_m
$,
 a probability measure concentrated on the subset $(A_z\times\{z\}) \cup (A_m\times\{m\})$
of~$A$. Notice that, for $i<n$,
$$
\nu\{y\in A:y_i\ne x_i\} = \th\mu_z\{w\in A_z: y_i\ne x_i\} +\thetabar\mu_m\{w\in A_m:
y_i\ne x_i\}
$$
 and
$$
\nu\{y\in A:y_n\ne x_n\} = 
\Cases
\thetabar&if $z\ne m$\\
0&otherwise
\endCases
\qquad\le\thetabar.
$$
By  the definition of~$F_n$ and the convexity of~$\psi$,
$$
\align
F_n(x,A) &\le \siln \psi\left(\strut \nu \{y_i\ne x_i\}\right)\\
&\le  \th\sum\NL_{i<n} \psi\left(\strut \mu_z \{y_i\ne x_i\}\right)
+\thetabar\sum\NL_{i<n} \psi\left(\strut \mu_m \{y_i\ne x_i\}\right)+\psi(\thetabar)
\\
\endalign
$$
The two sums over the first $n-1$ coordinates are like those that appear in the definitions
of~$F_{n-1}(w,A_z)$ and~$F_{n-1}(w,A_z)$. Indeed, taking an infimum over all 
$\mu_z\in\pp(A_z)$ and~$\mu_m\in\pp(A_m)$ we get the expression on the right-hand side
of~\cref{recursive.F}.

Take exponentials of both sides of~\cref{recursive.F} then integrate out \wrt/~$\QQ_{-n}$
over the~$w$ component. For $0<\th<1$ invoke 
the H\"older inquality, $\QQ_{-n} U^\th V^{\thetabar} \le  \left(\strut \QQ_{-n}
U\right)^\th
\left(\strut
\QQ_{-n} V\right){}^{\thetabar}$,
with $U=\exp( F_{n-1}(w,A_z))$ and $V=\exp( F_{n-1}(w,A_m))$, for a fixed~$m$.
For each~$z$ with~$A_z\ne\emptyset$ we get
\begin{equation}\label{QQn}
\begin{align}
&\QQ_{-n}\exp\left(\strut F_n((w,z),A)\right)\\
&\le \left(\strut \QQ_{-n}\exp\left(\strut F_{n-1}(w,A_z)\right)\right)^\th
\left(\strut \QQ_{-n}\exp\left(\strut F_{n-1}(w,A_m)\right)\right){}^{\thetabar}
e^{\psi(\thetabar)}\\
\end{align}
\end{equation}
The  inequality also hold in the extreme cases where $\th=0$ or $\th=1$,
by continuity. 
The inductive hypothesis bounds  
the last product by
$$
\exp\left(\strut \th\gam(R_z)+\thetabar\gam(R_m)+\psi(\thetabar)\right)
=\exp\left(\strut \gam(R_m) +G(\gam(R_z)-\gam(R_m),\th)\right)
$$
The exponent is a decreasing function of $R_m$.
Take an infimum over~$m$, to replace $\gam(R_m)$ by $\gam(r_0)$.
Then take an infimum over~$\th$ to get
\begin{equation}\label{QQn.bnd}
\align
\QQ_{-n}\exp&\left(\strut F_n((w,z),A)\right)\le 
\exp\left(\strut \gam(r_0)+G(\xi(R_z))\right)\\
&\text{where $\xi(r) := \gam(R_z)-\gam(r_0)$ for $0\le r\le
r_0$.}
\endalign
\end{equation}

If the crossection~$A_z$ is empty, the set~$\pp(A_z)$ is
empty. The argument leading from~\cref{recursive.F} to~\cref{QQn.bnd} still works if we
fix~$\th$ equal to zero throughout, giving the bound
$$
\QQ_{-n}^w\exp\left(\strut F_n(x,A)\right) \le \exp\left(\strut
\gam(r_0)+\psi(1)\right)
\qt{if $A_z=\emptyset$.}
$$
Thus the inequality~\cref{QQn.bnd} also holds with $R_z=0$ when~$A_z=\emptyset$,
because
$\xi(0)=\gam(0)-\gam(r_0)=\infty$ and $G(\infty)=\psi(1)$.

By Assumption~(i), the function $r\mapsto \exp\left(\strut
G(\xi(r))\right)$
 is  concave on~$[0,r_0]$. 
Integrate both sides of~\cref{QQn.bnd} with respect to~$Q_n$
to average out over the~$z$ variable.  Then invoke
Jensen's inequality and the fact that
$Q_n R_z =\QQ A$, to deduce that
$$
\QQ \exp\left(\strut F_n(x,A)\right)\le 
\exp\left(\strut \gam(r_0)+G\left(\strut \gam(\QQ A)-\gam(r_0) \right)\right).
$$
Finally, use the inequality $G(\eta)\le\eta$ to bound the last expression by
$\exp(\gam(\QQ A))$, thereby completing the inductive step.

\Remark 
Note that it is important to integrate \wrt/~$Q_n$ before using the bound on~$G$: the upper
bound $\exp(-\gam(R_z))$ is a convex function of~$R_z$, not concave.
\endRemark

\section{Proof of the Concavity Lemma}  \label{concave}
I will establish a more detailed set of results than asserted by
Lemma~\ref{concavity.lemma}.
Invoke the monotonicity and continuity of~$\psi'$ to define $g(\eta)$
as the solution to
$\psi'\left(\strut 1-g(\eta)\right) = \eta$ if $0\le\eta <\psi'(1)$
and~$g(\eta)=0$
if $\psi'(1)\le \eta$.
Then the following assertions are true.
%\end{document}
\begin{enumerate}
\ritem(i) 
$$
G(\eta) =\Cases
\psi\left(\strut 1-g(\eta)\right)+\eta g(\eta)& for $0\le\eta<\psi'(1)$\\
\psi(1)&for $\psi'(1)\le \eta$\\
\endCases  
$$

\ritem(ii) $G$ 
is increasing and concave, with a continuous, decreasing  first derivative~$g$. 
In particular,  $G(0)=0$ and
$G'(0)=g(0)=1$.

\ritem(iii)
$
G''(\eta) =g'(\eta)= -\left[\psi''\left(\strut 1-g(\eta)\right)\right]^{-1}$
for $0<\eta<\psi'(1)$.

\ritem(iv)  $G(\eta)\le \eta$ for all $\eta\in\RR^+$.

\ritem(v)
Suppose $\xi:J\to\RR^+$ is a convex function defined on a subinterval~$J$ of the real line,
with $\xi'\ne0$  on the interior of~$J$. 
Suppose
$$
{1\over \psi''(1-\xi_r)} \ge g(\xi_r)^2 + g(\xi_r) \xi''(r)/\xi'(r)^2,
$$ 
for all $r$ in the interior of~$J$ for which $\xi_r :=\xi(r)\in(0,1)$. 
Then 
$r\mapsto \exp\left(\strut G(\xi(r))\right)$ is a concave function on~$J$.
\end{enumerate}

\demo{Proof of (i) through (iv)} 
The fact that $G$ is concave and increasing follows from its definition as an infimum of
increasing linear functions of~$\eta$. (It would also follow from the fact that
$G'(\eta)=g(\eta)$, which is nonnegative and decreasing.) Replacement of the infimum over
$0\le \th\le1$ by the value at~$\th=1$ gives the inequality~$G(\eta)\le
\eta$.

If $\eta\ge \psi'(1)$, the derivative $-\psi'(1-\th)+\eta$ is nonnegative on~$[0,1]$,  which
ensures that the infimum is achieved at $\th=1$. 

If $0<\eta<\psi'(1)$, the infimum is achieved at the zero of the derivative, $\th=g(\eta)$.
Differentiation of the defining equality $\psi'\left(\strut 1-g(\eta)\right) = \eta$ then
gives the expression for~$g'(\eta)$. Similarly
$$
G'(\eta) = -\psi'\left(\strut 1-g(\eta)\right) g'(\eta) +\eta g'(\eta) + g(\eta) = g(\eta).
$$

The infimum that defines~$G(0)$ is achieved at $g(0)=1$, which gives $G(0)=\psi(0)=0$.
Continuity of~$g$ at then gives $G'(0)=g(0)=1$.

\demo{Proof of (v)}  Note that the function $L(r) := \exp\left(\strut G(\xi(r))\right)$
is continuous on~$J$ and takes the value $e^{\psi(1)}$ for all~$r$ at which $\xi(r)\ge
\psi'(1)$. The second derivative~$L''(r)$ exists except possibly at points~$r$ for
which~$\xi(r)=\psi'(1)$.  In particular, $L''(r)=0$ when $\xi(r)>\psi'(1)$ and
$$
L''(r) = \left(\strut g'(\xi_r) (\xi'_r)^2 +g(\xi_r)\xi''_r +
g(\xi_r)^2(\xi'_r)^2\right)L(r)
\qt{for }0<\xi_r<\psi'(1).
$$
From~(iii) and the positivity of~$L$, the last expression is $\le0$ if and only if
$$
-{(\xi'_r)^2\over\psi''(1-g(\xi_r))}  +g(\xi_r)\xi''_r +
g(\xi_r)^2(\xi'_r)^2 \le 0
$$
Divide through by~$(\xi'_r)^2$ then rearrange to get the asserted inequality for $\psi''$.
%\endProof
Lemma~\ref{concavity.lemma} follows as a special case of~(i) through~(iv).

%\bigskip
%\hbox to\hsize{\hfil\qed\qquad \qed\qquad \qed\hfil}
%\bigskip

\bigskip

\demo{Special cases}
If $\sup_r\xi''(r)/\xi'(r)^2\le c$, with~$c$ a positive constant,
 the
inequality from part~(v) will certainly hold if
\begin{equation}
\psi''(1-\th) \le (\th^2 + c\th )^{-1}
\qt{for all }0<\th<1.\label{diff.ineq}
\end{equation}
This differential inequality can be solved, subject to the constraints
$0=\psi(0)=\psi'(0)$, by two integrations.  
Then
$$
\psi'(1-\th) =\int_\th^1 \psi''(1-t)\,dt
 \le \int_\th^1 \frac{dt}{t^2+ct} 
= c^{-1}\left(-\log \th + \log\left(\strut {\th+c\over 1+c}\right)\right)\\
$$
and, with $\psi_c$ defined by~\cref{psic},
$$
\psi(1-\th) = \int_\th^1 \psi'(1-t)\,dt 
\le c^{-1}\int_\th^1 -\log t + \log\left(\strut {t+c\over 1+c}\right)\,dt
=\psi_c(1-\th)  .
$$
Note that $\psi_c(1-\th)$ is the solution to the differential
equation 
$$
\psi_c''(1-\th) = {1\over \th^2 + c\th }
\qt{for all $0<\th<1$, with $\psi_c(0)=\psi_c'(0)=0$.}
$$
It is the largest solution to~\cref{diff.ineq}.

%\bibliographystyle{chicago}
%\bibliography{dbp}

\def\noopsort#1{}

\end{document}